\documentstyle[12pt]{article}
\begin{document}
\def\square#1{\vbox{\hrule
\hbox{\vrule\hbox to #1 pt{\hfill}\vbox{\vskip #1 pt}\vrule}\hrule}}
\newtheorem{guess}{Proposition }[section]
\newtheorem {theorem}[guess]{Theorem}
\newtheorem{lemma}[guess]{Lemma}
\newtheorem{corollary}[guess]{Corollary}
\newtheorem{example}[guess]{Example}
\newtheorem{remark}[guess]{Remark}
\def \ja {\vrule height 3mm width 3mm}

\centerline{ \Large \bf Numerical invariants for bundles on blow-ups}
\vspace{.1 in}
\centerline{\sc E. Ballico{\footnote{Partially supported by MURST (Italy)}} - E.   Gasparim{\footnote{Partially supported
by CNPQ (Brasil)}}}

\vspace{.2in}


\begin{abstract}

We suggest an effective procedure to 
calculate numerical invariants for  rank two  bundles
over blown-up surfaces. 
We study 
the moduli spaces ${\cal M}_j$  of rank two  bundles 
on the blown-up plane splitting over the exceptional divisor
as ${\cal O}(j) \oplus {\cal O}(-j).$
We use the 
 numerical invariants
to give a topological decomposition of ${\cal M}_j.$

\end{abstract}

\section{ Introduction}

Let $X$ be a complex surface and let $\pi : \widetilde {X} \rightarrow X$
be the blow-up of $X$  at point $p \in X.$ We denote by $\ell$ the 
exceptional divisor.
For simplicity we  restrict ourselves to 
bundles with vanishing first Chern class 
(although our methods  apply 
to more general situations).

If $E$ is a rank two 
bundle over a surface $X$
and $j$ is an integer, then a  polynomial 
$p$ determines a bundle $\widetilde{E}$ 
over the blown-up surface $\widetilde{X}$
such that $\pi_*\widetilde{E}^{\vee\vee} = E$
and $\widetilde{E}|_{\ell} = {\cal O}(j) \oplus {\cal O}(-j).$
It follows that we have a correspondence 
 $\{p\} \mapsto \{rank\,\, two\,\,bundle\,\,
\widetilde{E} \rightarrow \widetilde{X}\}$(see \cite{BSPM}).
 We 
calculate numerical invariants for 
$\widetilde{E} $
as a function of the polynomial $p$
 and in particular we
 give an effective procedure
to calculate the charge   $c_2(\widetilde{E}) -
c_2(E)$
which is the difference of second  Chern classes.
This charge
is a sum of two finer  invariants  $ l(Q)$ and
$ l(R^1\pi_*\widetilde{E}),$ 
where $Q$ is the sheaf defined by the exact sequence 
$0 \rightarrow \pi_*\widetilde{E} \rightarrow 
( \pi_*\widetilde{E})^{\vee\vee}
 \rightarrow Q \rightarrow 0$
(see \cite{FM}).

Consider the moduli spaces ${\cal M}_j$ of bundles 
on the blown-up plane which split over 
the exceptional divisor as ${\cal O}(j) \oplus {\cal O}(-j).$
There is a natural quotient topology on ${\cal M}_j$  
(see  \cite{JA}) 
which turns out to be non-Hausdorff. This is to be expected since 
$\widetilde{\bf C}^2$ is non-compact and 
there is a priori no notion of stability for bundles 
over a non-compact space. However, we show that the numerical 
invariants determine a nice decomposition of  ${\cal M}_j$
into Hausdorff subspaces.

\section{A correspondence between polynomials and vector bundles}

Consider a rank two holomorphic bundle $\widetilde{E}$ on $\widetilde{X}$ 
with vanishing first 
Chern class.
Let  $j \geq 0$   be the splitting type  of $\widetilde{E}$ to
the exceptional divisor, i.e.,
$\widetilde{E}|_{\ell} = {\cal O}(j) \oplus {\cal O} (-j).$
  Then on a neighborhood $N$ of the
exceptional divisor,
there is a canonical choice of transition 
matrix for $\widetilde{E}|_N,$
 namely
$T = \left(\matrix {z^j & p \cr 0 &  z^{-j} \cr }\right),$
where 
$$(*) \hspace{3cm} p = \sum_{i =1}^{2j -2}\sum_{l = i-j+1}^{j-1} p_{il}z^lu^i $$
is a polynomial with  complex coefficients 
(see \cite{JA}).
In other words, $\widetilde{E}|_N$ is 
given as an algebraic extension
$0 \rightarrow {\cal O}(-j) 
\rightarrow \widetilde{E}|_N \rightarrow
{\cal O}(j) \rightarrow 0$
 whose extension class is determined by $p,$
or equivalently we may say that  $\widetilde{E}|_N$
is determined by the pair $(j,p).$
It follows that every 
holomorphic rank two vector bundle over $\widetilde{X}$
with vanishing first Chern class 
is topologically  determined 
 by a triple $(E,j,p)$ (see \cite{BSPM}). 

If we fix an integer $j$ 
and a holomorphic bundle $E$ over 
$X,$ then we get a correspondence 
between polynomials and bundles
on $\widetilde{X}.$
$$
\begin{array}{rcl}  & \phi & \cr
\left\{ polynomials\,\, of\,\, the\,\, form\,\, (*)\right\} &  \rightarrow &
\left\{ rank \,\,two \,\, bundles \,\, on \,\, \widetilde{X}\right\} \cr
 p  & \mapsto &  \left(\matrix {z^j & p \cr 0 &  z^{-j} \cr }\right).
\end{array}$$
Here $im(\phi)$ consists of all bundles $\widetilde{E}$ over 
$\widetilde{X}$ 
satisfying the following two conditions:

$\iota$) $\widetilde{E}|_{\ell} = {\cal O}(j) \oplus {\cal O}(-j)$

$\iota\iota$) $\widetilde{E}|_{\widetilde{X} - \ell} = \pi^* (E|_{X-p}).$
 
Note that varying $E$ we get all holomorphic bundles 
over $\widetilde{X}$ which split as ${\cal O}(j) \oplus {\cal O}(-j)$
over the exceptional divisor and  if we also vary $j$ 
we obtain all rank two holomorphic bundles over $\widetilde{X}$ with vanishing
first Chern class (see \cite {BSPM}). 

\section{Numerical Invariants}

Let $\widetilde{E}$ be any holomorphic
rank two bundle over 
$\widetilde{X}$ with vanishing first Chern class.
In a neighborhood of the exceptional divisor, the bundle      
 $\widetilde{E}$ has a canonical transition 
 matrix $T,$ which in turn defines a unique bundle 
$\widetilde{V}$ over the blown-up plane.
On the other hand, if we start with a bundle $\widetilde{V}$ 
on the blown-up plane and a bundle $E$ over $S,$ we 
show in section II that we can 
glue them together to form a bundle $\widetilde{E}$ over 
$\widetilde{S}.$ 
As mentioned in the introduction, the difference of second 
Chern classes is given by two numerical invariants:
 $c_2(\widetilde{E}) - 
c_2(E)= l(Q) + l(R^1\pi_*\widetilde{E}).$ 
But these are local invariants and therefore  can be viewed as
invariants corresponding to the bundle $\widetilde{V}$
itself. 

Recall that we say that $\widetilde{E}$ has splitting type 
$j$ when $\widetilde{E}|_{\ell}
\simeq {\cal O}(j)\oplus {\cal O}(-j).$
For splitting types $0$ and $1$ the
moduli spaces ${\cal M}_0$ and ${\cal M}_1$ are trivial.
In fact,
every bundle over $\widetilde{\bf C}^2$ which is trivial over
the exceptional divisor is trivial over the entire $\widetilde{\bf C}^2,$
hence ${\cal M}_0$ is just a point; and every bundle on
over $\widetilde{\bf C}^2$ with splitting type $1$
over $\ell$ splits over the entire 
$\widetilde{\bf C}^2,$  therefore ${\cal M}_1$ 
is also just a point (see \cite{JA}). Therefore the numerical invariants 
are also trivial in these two cases: both invariants vanish for
splitting type zero and for splitting type one we have
$l(R^1\pi_*\widetilde{E})= 0$ and $l(Q)= 1.$ 

Numerical invariants for splitting types 2 and 3 
are tabulated below. We leave the calculations for the last section.
The notation for the tables is the following.
For fixed splitting type,
each bundle $\widetilde{V}$
is given as an extension of 
line bundles with the extension class determined 
by the  polynomial $p$ written in the canonical form $(*)$ 
as in section 2. We write the coefficients of 
$p$ in lexicographical order of the coefficients $p_{il}$ and 
 denote by $e_n$ the monomial corresponding to the $n-th$ term.
For instance for splitting type $2$ the polynomial written in this order is 
$p = p_{10}u+p_{11}zu+ p_{21}zu^2$ and $e_1 = u,$
$e_2= zu,$ $e_3= zu^2.$
We use the term {\em charge} for the difference 
$c_2(\widetilde{E}) - c_2(E).$ 

\vspace{5 mm}

\begin{tabular}{|c|c|c|c|}  \hline
TABLE  I &\multicolumn{3}{|c|}{$splitting\,\, type\,\,  2$}\\
\hline
{\em monomial} & {\em l(Q)} & $ l(R^1\pi_* \widetilde{V})$ & charge \\
\hline
 $e_1$ & 1  &  1 & 2 \\
  $e_2$ & 1 & 1 & 2 \\
\hline
  $e_3$  & 2 & 1 & 3 \\
\hline
 zero & 3 &  1 & 4\\
 \hline
\end{tabular}

\vspace{5 mm}

\begin{tabular}{|c|c|c|c|}  \hline
TABLE II &\multicolumn{3}{|c|}{$splitting \,\, type\,\,  3$}\\
\hline
{\em monomial} &{\em l(Q)}& $l(R^1\pi_*\widetilde{V})$ & charge \\
\hline
 $e_1$ &  3 & 2 & 5 \\
  $e_2$ & 1 & 2 & 3 \\
  $e_3$ & 1 & 2 & 3 \\
  $e_4$ & 3 & 2 & 5 \\
\hline
  $e_5$ & 3 &  3 & 6\\
  $e_6$ & 2 &  3 & 5 \\
  $e_7$ & 3 &  3 & 6 \\
\hline
  $e_8$ & 4 &  3 & 7\\
  $e_9$ & 4 &  3 & 7\\
\hline
  $e_{10}$ & 5 & 3 &  8\\
\hline
  zero &  6 & 3 & 9\\
 \hline
\end{tabular}

\vspace{5 mm}

\begin{tabular}{|c|c|c|c|}  \hline
TABLE  III &\multicolumn{3}{|c|}{$splitting\,\, type\,\,  3$}\\
\hline
{\em polynomial} & {\em l(Q)} & $ l(R^1\pi_* \widetilde{V})$ & charge \\
\hline
 $e_1+e_4$ & 1  &  2 & 3 \\
 \hline
  $e_4+e_5$  & 2 & 2 & 4 \\
\hline
 $e_1+e_7$ & 2 &  2 & 4\\
 \hline
\end{tabular}

\section{Moduli of  bundles on the blown-up plane}

Let ${\cal M}_j$ denote the moduli space of  
 bundles 
over the blown-up plane whose restriction to the exceptional divisor equals
$ {\cal O}(j) \oplus {\cal O}(-j).$
Our goal here is to show that the two  numerical invariants 
from the previous section provide a decomposition
 of ${\cal M}_j$ into Hausdorff subsets.

A remark about the terminology
used here is in order. Rigorously we should
use the term {\em parameter spaces}
for the spaces ${\cal M}_j,$
since among other things these spaces are non-Hausdorff
and since there is no notion of stable points in this case.
One possible choice for a notion of stability would be to
call a bundle stable when it  belongs to the generic subset
of ${\cal M}_j,$ which in terms of the numerical invariants corresponds
to having the smallest possible invariants.
Since this generic part is  Hausdorff, these
would seem appropriate, but it is not yet clear
if this is the best notion of stability to impose.

First we recall the topology of ${\cal M}_j.$
Each
 bundle $E$ on $\widetilde{\bf C}^2$ with splitting type $j$ 
is represented  by a pair $(j,p)$ where $p$ 
determines the extension class (see section 2).
Writing the polynomial  $\sum p_{il}z^lu^i$ in lexicographical order, 
gives a natural  
identification of $p$  with the point in ${\bf C}^N$
whose coordinates are the  coefficients of $p.$
It is natural to  we impose the relation $p \sim p'$ if 
$(j,p)$ and $(j,p')$ represent  isomorphic bundles, and 
take  ${\bf C}^N/ \sim$
with the quotient topology.
We give ${\cal M}_j$ the topology induced by the bijection
${\cal M}_j \rightarrow {\bf C}^N/ \sim.$

It turns out that once this quotient  topology is described, 
a  decomposition of the moduli spaces 
essentially suggests itself, and nicely enough 
this decomposition is the same as the one we obtain by
separating loci of constant pairs of numerical invariants
$(l(Q),l(R^1\pi_*\widetilde{V})).$
In general this decomposition becomes a bit 
abstract, but it is quite clear in the  first  examples
where we take the splitting type to be a small integer.

\subsection {The topological  structure of ${\cal M}_2$}

As mentioned in section 3, the
 moduli spaces ${\cal M}_0$ and ${\cal M}_1$
are trivial; both consist of a single point. The first interesting
example happens when the splitting type is
$2,$ where we have the following structure
(see \cite{JA}).

The moduli space ${\cal M}_2$ is isomorphic to 
${\bf C}^3/\sim$ where $\sim$ is
the equivalence relation given by  $R_1,$ and  $R_2$ 

\begin{itemize}
\item $R_1:$ $(p_{10},p_{11},p_{21}) \sim  \lambda (p_{10}, p_{11},
*)$ if $(p_{10},p_{11}) \neq (0,0)$ 
\item $R_2:$ $(0,0,p_{12}) \sim  \lambda ( 0,0,  p_{12})$
\end{itemize}
where $ \lambda \in {\bf C} -\{0\}.$
This   
suggests a decomposition 
${\cal M}_2 = S_0 \cup S_1 \cup S_2$
where 
\begin{itemize}
\item $S_0= \{(0,0,0)\}$
\item $S_1  =\{(\lambda p_{10},\lambda p_{11}, *),
 (p_{10},p_{11}) \neq (0,0)\} 
 \simeq {\bf P}^1 $ 
\item   $S_2  = \{\lambda (0,0,p_{21}),p_{21} \neq 0\} \simeq \{(0,0,1)\}$
  \end{itemize}

On the other hand, one could also use the numerical invariants 
given on Table I to decompose ${\cal M}_2.$
Comparing  with Table I 
we see that the  $S_i \subset {\cal M}_2$ with $i = 0,1,2$
coincide with the loci of constant $l(Q)$ (or else the loci of 
constant charge).

\subsection {The topological  structure of ${\cal M}_3$}

For splitting type $3$ the quotient space becomes 
somewhat more complicated, and  both numerical invariants
are needed to give a nice decomposition.
In the special case of ${\cal M}_2$ 
we saw that the decomposition
by loci of constant charge already provided 
us with Hausdorff subspaces $S_i.$
However, this is in general not the case and as we 
shall see in the example of ${\cal M}_3.$
In fact, the two finer invariants $l(Q)$ and $l(R^1\pi_*\widetilde{V})$
are needed two give a nice decomposition;
or equivalently, one can choose
 to give the charge and one of the finer invariants.

 For splitting type $3$
 moduli space is  ${\cal M}_3 \simeq {\bf C} ^{10} / \sim $
with equivalence relation defined by  $R_i, i = 1,...,6$ below,
where $\lambda \in {\bf C}-\{0\}$ and  $*$  denotes either 
a complex   number or a few  complex numbers,
on which no restrictions are imposed

\vspace{5mm}
\noindent $\bullet \,\, R_1:$ $(a_1,\cdots,a_{10})  
  \sim \lambda ( a_1 , \cdots,
a_4, *)$ if $(a_2,a_3) \neq (0,0)$ or  $a_1,a_4 \neq 0$
 
\vspace{3mm}
\noindent $\bullet \,\, R_2:$ $(a_1,0,0,0,a_5,\cdots,a_{10})  \sim
\lambda(a_1 ,0,0,0, *, a_7, *)$ 
if $a_1 \neq 0$\\
\hspace*{1.2cm}$(0,0,0,a_4,\cdots,a_{10})  \sim
\lambda(0,0,0 ,
a_4, a_5, * )$
if $a_4 \neq 0$

\vspace{3mm}
\noindent $\bullet \,\, R_3:$
$  ( 0,\cdots,0,a_5,\cdots,a_{10})  \sim
\lambda(0,\cdots,0 ,a_5,a_6,
a_7, * )$ 
if $a_6 \neq 0$ or  $a_5,a_7 \neq 0$

\vspace{3mm}
\noindent $\bullet \,\, R_4:$
$  (0,\cdots,0,a_5,\cdots,a_{10})  \sim
\lambda(0,\cdots,0,a_5,0,
0 ,*,a_9,*)$ 
if $a_5 \neq 0$\\
\hspace*{1.2cm}$  (0,\cdots,0,a_7,\cdots,a_{10})  \sim
\lambda(0,\cdots,0,
a_7 ,a_8,*)$ 
if $a_7 \neq 0$

\vspace{3mm}
\noindent $\bullet \,\, R_5:$
$  (0,\cdots,0,a_8,a_9,a_{10})  \sim
\lambda(0,\cdots,0,a_8,a_9, *)$ 
 if  $ (a_8,a_9) \neq  (0,0) $

\vspace{3mm}
\noindent $\bullet \,\, R_6:$
$  (0,\cdots,0,a_{10})  \sim
\lambda(0,\cdots,0,a_{10})$  

\vspace{5 mm}

This equivalence relation suggests a decomposition of the space ${\cal M}_3$ 
into subsets $S_i.$ Let us look at the 
most generic part of ${\cal M}_3.$
>>>From the first relation $R_1$ it is natural to consider 
 the subset $X_1 =  {\bf C}^{10} -  V({\cal I})$
where ${\cal I} = <x_2,x_3,x_1x_4>.$ Then the 
subspace $S_1 = X_1/R_1$ is homeomorphic to 
${\bf CP}^3- \{[1,0,0,0], [0,0,0,1]\}.$ 
 Comparing with tables II and III we see 
that for points in $S_1$
the lowest values of the numerical invariants occur.
In fact, $S_1$ is exactly the locus of points in ${\cal M}_3$  
attaining  the lowest  invariants.

On the opposite side, relation $R_6$ 
gives us the point 
$  (0,\cdots,0,1) $ which is the only point 
of ${\cal M}_3$  corresponding to charge 8.
There is one  point which is more special.  The single  point  
$\{0\}$ forms the 
least generic subset of ${\cal M}_3$ 
 where  the highest values
 of the numerical invariants are attained.

Among the other relations we notice 
another interesting fact. Let us look at relations $R_2$ and $R_3.$
>>>From relation $R_2$ there comes out a subset $S_2$  of ${\cal M}_3$
  homeomorphic
to two copies of ${\bf CP}^1 - [0,1].$
On the other hand, from relation $R_3$ there comes out a subset $S_3$
of ${\cal M}_3$
homeomorphic to ${\bf CP}^2 - \{[1,0,0],[0,0,1]\}.$
Naturally we would like to have $S_3$ and $S_2$ as separate subsets on 
a nice decomposition of ${\cal M}_3,$ specially because 
 their union is non-Hausdorff in ${\cal M}_3.$
However, if we were to decompose ${\cal M}_3$ by loci of constant charge, then 
$S_2 \cup S_3$ would be contained in the single subset of 
points corresponding to charge 5. This makes it evident that 
one more numerical invariant is necessary to give a nice decomposition
of ${\cal M}_3,$ and it is clear that any pair of  invariants 
from table II  distinguishes $S_2$ from $S_3.$

It comes out nicely, that the topological decomposition 
that is naturally suggested by the description of 
${\cal M}_3$  as a quotient of ${\bf C}^{10}$ 
is exactly the same as the decomposition given by loci of 
constant pairs of numerical invariants
$(l(Q),l(R^1\pi_*\widetilde{V})).$

\subsection{The topological structure of ${\cal M}_j.$}

The facts mentioned about the topology of ${\cal M}_3$ are readily 
generalized for higher splitting type.
The key fact to have in  mind is  that
to calculate the numerical invariants, one takes into account 
which are the nonzero coefficients of the polynomial, but not the particular
value of each coefficient. 

The generic set of ${\cal M}_j$ 
is a homeomorphic to a 
 complex projective space of dimension $2j-3$ minus a complex subvariety 
of codimension at least 2 (see \cite{JA}). For points in this generic part,
the lower bound of the numerical invariants are attained 
and these are (see \cite{TO}) $l(Q) = 1$ and 
$l(R^1\pi_*\widetilde{V}) = j-1$ hence charge $j.$ 
The least generic point in ${\cal M}_j$
comes from $0 \in {\bf C}^N$ 
which by the correspondence from section 2 
gets translated into the vanishing of the 
 extension class,  that is, to the split bundle
 ${\cal O}(j) \oplus {\cal O}(-j),$
for which we have (see \cite{TO})
$l(R^1\pi_* \widetilde{V}) = j(j+1)/2$
and $l(Q) = j(j-1)/2$ and hence charge $j^2.$
Every intermediate value of the numerical invariants occur for
some point in ${\cal M}_j$ (see \cite{BG}).
To show that these invariants provide a nice decomposition of 
${\cal M}_j$ into Hausdorff subspaces, 
we use induction over $j$ together with the fact that 
there is a topological embedding ${\cal M}_{j-1} \hookrightarrow {\cal M}_j$
(see \cite{HO}).
 We then have a topological decomposition of ${\cal M}_j$
into a union of subspaces homeomorphic to
  open subsets of  complex projective spaces
${\bf CP}^n$ with $0\leq n\leq 2j-3$  and two points.
In other words, we have just showed:

\begin{theorem} The numerical invariants $(l(Q),
l(R^1\pi_* \widetilde{V}))$ provide a  decomposition ${\cal M}_j
= \cup S_i$ where each $S_i$ is homeomorphic to an open subset 
of a complex projective space of dimension at most $2j-3.$
 The lower bounds for these invariants 
are $(1,j-1)$ and this pair of invariants takes
  place on the generic part of ${\cal M}_j$ 
which is homeomorphic to   ${\bf CP}^{2j-3}$ 
minus a closed subvariety of 
codimension at least 2. 
The   upper bounds for these invariants are $(j(j-1)/2, j(j+1)/2)$
and this pair occurs at one single point of ${\cal M}_j$ 
which represents the split bundle.
\end{theorem}

\section{Computing  the invariants}
The aim of this section is to  give a more concrete feeling
about the invariants, their  geometric meaning and how to 
calculate them. We remark that the calculations follow an
algebraic procedure which has as initial data just the
transition matrix $T$ for the bundle $\widetilde{V}$ 
on a neighborhood of the exceptional divisor. Given this data 
one calculates the zero-th and first cohomology groups
for $\widetilde{V}$ 
and then 
it is a matter of simple algebra to calculate 
$l(Q)$ and $l(R^1\pi_*\widetilde{V}).$ 
Despite the fact that the calculations 
are long, they are quite simple. In principle 
one could write a computer program 
to do them, and it would be 
 interesting to 
tabulate the invariants for higher values of 
the splitting type $j$
as it would give a better feeling for how the 
Hausdorff  subsets
 of ${\cal M}_j$
are distributed. 
Since the sum of these numbers 
gives the charge,
in particular this method gives a completely
algebraic procedure to calculate the second Chern class of $\widetilde{V}$
directly from its transition matrix.
Calculations for first Chern class of line bundles 
from the transition matrix are well known,
but for rank two bundles the authors do not
know of any reference.

\subsection{Geometric meaning}
Let us  make some comments about the 
geometric interpretation of these numbers.
Recall that these are local invariants, so that 
the geometric meaning corresponds to 
the behavior of the bundle $\widetilde{E}$
in a neighborhood of the blow-up.
Suppose we are given the transition matrix $T= 
\left(\matrix{z^j & p \cr 0 & z^{-j}}\right)$ for
$\widetilde{V}.$
The data on this matrix means that
$\widetilde{V}$ is given as an extension of line bundles
$0 \rightarrow {\cal O}(-j) \rightarrow \widetilde{V} \rightarrow {\cal O}(j)
\rightarrow 0$ with extension class determined by the 
polynomial $p.$
When $p = 0$ the bundle  splits and  
$l(R^1\pi_*\widetilde{V}) = j(j-1)/2$ assumes its maximal value.
  For the most general cases $p$
is nonzero on the first formal neighborhood
and  $\widetilde{V}$ 
belongs to the generic part of 
${\cal M}_j$
in which case $l(R^1\pi_*\widetilde{V}) = j-1$
takes the lowest value.
The difference $j(j-1)/2 - l(R^1\pi_*\widetilde{V})$ is a measure of 
``how far'' the bundle is from 
being split.

The invariant $l(Q)$ is the length
of the sheaf $Q$ defined by the exact sequence
$0 \rightarrow  \pi_*\widetilde{E} 
\rightarrow (\pi_*\widetilde{E})^{\vee\vee}
\rightarrow Q \rightarrow 0.$
Note that $Q$ is supported only at a point.
If $ \pi_*\widetilde{E}$ is locally free then $Q$
is trivial and $l(Q) = 0,$  but this only happens 
when $j = 0$ in which case $\widetilde{E}$ is a pull back;
otherwise,  $ \pi_*\widetilde{E}$ is not locally free.
The length 
$l(Q)$ is the dimension of the stalk at this point
and it measures ``how far'' the sheaf  $ \pi_*\widetilde{E}$
is from being locally free, which can also be seen as a
 measure of  ``how far'' $\widetilde{E}$ is from being a pull back bundle.

\subsection{ How to calculate  $l(Q)$}

Let $M = (\pi_* \widetilde{V})_x^{\wedge}$ denote the
completion of the  stalk 
 $(\pi_*\widetilde{V})_x$ over the blown-up point $x.$
Let $\rho$ denote the natural inclusion of $M$ 
into its bidual $\rho: M \hookrightarrow M^{\vee\vee}.$
We want to compute $l(Q) = dim\, coker(\rho).$ 
By the theorem on formal functions (see \cite{HA})
$$M \simeq   \lim_{\longleftarrow}H^0(\ell_n, \widetilde{V}|{\ell_n}).$$
There are  simplifications that make it easy to  calculate $M.$
For a fixed splitting type  $j$
 it is sufficient to calculate 
$H^0(\ell_n, \widetilde{V}|{\ell_n})$ for $n 
\leq 2j-2.$
This follows from the fact that the polynomial $p$ 
determining the extension class
 has nonzero coefficients  only up to the
$(2j-2)-{nd}$ formal neighborhood. 
Moreover, the groups  $H^0(\ell_n, \widetilde{V}|{\ell_n})$
and $ H^0(\ell_{n+1},\widetilde{V}|{\ell_{n+1}})$
 for $n>2j-2$
have the same generators as ${\cal O}_x$-modules.
 It follows that to determine  $M$ 
 it suffices to
 calculate  $H^0(\ell_{2j-2}, \widetilde{V}|{\ell_{2j-2}}),$
 and the relations among its generators
under the action of   ${\cal O}^{\wedge}_x ( \simeq
{\bf C}[[x,y]]$).
In what follows we fix a coordinate system for
$\widetilde{\bf C}^2$ given by two charts 
$U = \{(z,u)\}\simeq {\bf C}^2 \simeq V = \{(\xi,v)\}$
with $(\xi,v) = (z^{-1},zu)$ on $U\cap V.$
Since the blow-up map
$\pi : \widetilde{{\bf C}^2} \rightarrow {\bf C}^2$
is given  by 
$ (x,y) = \pi(z,u) = (u,zu)$ on the $U$ chart
the  natural action of $x$  and $y$ on this space is that  $x$
acts by multiplication by $u$ and 
$y$ acts by multiplication by $zu.$

Calculations  of $l(Q)$ for the split case and also for the 
generic case (corresponding to $p = u$) appeared in \cite{TO}.
Here we present  the calculations
for the bundle given by transition matrix
$T = \left( \matrix{ z^2 & zu^2 \cr
                       0 & z^{-2} } \right)$
which is neither generic nor split and corresponds to the 
subset $S_2 \subset {\cal M}_2.$
According to our previous remarks it suffices to
 calculate the generators 
of  $H^0(\ell_2, \widetilde{V}|{\ell_2})$
and the relations among them, from which we find that

$M = {\bf C}[[x,y]]<\alpha_0,\beta_0,\beta_1,\beta_2>$
where 

$$
\alpha_0 = \left(\matrix{ u^2 \cr 0}\right),\,\,\,
\beta_0 = \left(\matrix{ 0 \cr 1}\right),\,\,
\beta_1 = \left(\matrix{ 0 \cr z}\right),\,\, 
\beta_2 = \left(\matrix{ -zu^2 \cr z^2 }	\right)
$$
with relations
$$\left\{\matrix {x\beta_1 - y \beta_0 \cr
x \, \beta_2 - y \, ( \alpha_0 +\beta_1)  }\right..$$
Once this is found, it is a simple algebraic calculation to 
find the dual and bidual.
We have that 
$M^{\vee} = < A ,  B , C >$  has the generators 
$$  A  :\left\{\matrix {\alpha_0 \rightarrow x \cr
\beta_2 \rightarrow y \cr
\beta_0 \rightarrow 0 \cr
\beta_1 \rightarrow 0 \cr
}\right. \,\,\, 
 B :\left\{\matrix {\alpha \rightarrow 0 \cr
\beta_2 \rightarrow y^2 \cr
\beta_0 \rightarrow x^2\cr
\beta_1 \rightarrow xy \cr
}\right. \,\,\,
 C :\left\{\matrix {\alpha \rightarrow -y \cr
\beta_2 \rightarrow 0 \cr
\beta_0 \rightarrow x \cr
\beta_1 \rightarrow y \cr
}\right. 
 $$
satisfying the relation
$$ y \, A - B + x \, C = 0.$$
And   $M^{\vee\vee}$  is free on two generators
(this will be always the case since it is 
the stalk of a rank two locally free sheaf).
$M^{\vee \vee} = <{\cal A}, {\cal B}>$
where
$$ {\cal  A}  :\left\{\matrix { A \rightarrow 1 \cr
B \rightarrow y \cr
C \rightarrow 0 \cr
}\right. \,\,\, 
 {\cal B} :\left\{\matrix {A \rightarrow 0 \cr
B \rightarrow x \cr
C \rightarrow 1 \cr
}\right..$$

The map $\rho :M \rightarrow M^{\vee \vee},$
 is given by evaluation.
We have
$$ \rho :\left\{\begin{array}{l} 
\alpha_0 \rightarrow  x{\cal A} - y{\cal B}  \cr
\beta_2 \rightarrow y {\cal A} \cr
\beta_0 \rightarrow x {\cal B} \cr
\beta_1 \rightarrow y {\cal B}
\end{array}\right. $$
Hence $ im (\rho) = < x {\cal A} - y {\cal B},
y {\cal A} , x  {\cal B} , y  {\cal B} >,$ and  $coker( \rho) = 
< \overline{\cal A} , \overline{\cal B} >$
and therefore $l(Q) = dim\, coker( \rho) = 2.$

\subsection{How to calculate $l(R^1\pi_*\widetilde{V})$}

Using the theorem on formal functions we have 
 $$l(R^1\pi_*\widetilde{V}) = 
 dim\, \lim_{\longleftarrow}H^1(\ell_n, \widetilde{V}|{\ell_n}).$$
But because the extension class for $\widetilde{V}$ is
given by the polynomial $p$ which has nonzero coefficients only
up to the $(2j-2)-{nd}$ formal neighborhood it suffices to calculate 
$H^1(\ell_{2j-j}, \widetilde{V}|{\ell_{2j-2}}).$  
The numerical invariant $l(R^1\pi_*\widetilde{V})$ counts just the 
number of generators 
of $H^1(\ell_{2j-2}, \widetilde{V}|{\ell_{2j-2}}).$ 
The calculations of  $l(R^1\pi_*\widetilde{V})$ are even simpler than 
those for $l(Q).$

We take the same example as in the previous section.
Let $\widetilde{V}$ be given by transition matrix
$T = \left( \matrix{ z^2 & zu^2 \cr
                       0 & z^{-2} } \right).$
We need to find the generators for  $H^1(\ell_2, \widetilde{V}|{\ell_2}).$
Let $\sigma \in  H^1(\ell_2, \widetilde{V}|{\ell_2}) $ 
then $\sigma = \sum_{i = 0}^2\sum_{z= -\infty}^{\infty}
\left(\matrix{a_{ik} \cr b_{ik}}\right)z^ku^i.$
But $ \sum_{i = 0}^2\sum_{z= 0     }^{\infty}
\left(\matrix{a_{ik} \cr b_{ik}}\right)z^ku^i $
gives a holomorphic function in $U$ and therefore represents a
 coboundary, which may subtract from $\sigma$
without changing its cohomology class.
Hence $\sigma \sim   \sum_{i = 0}^2\sum_{z= -\infty}^{-1}
\left(\matrix{a_{ik} \cr b_{ik}}\right)z^ku^i.$
Changing coordinates we have $T \sigma = 
   \sum_{i = 0}^2\sum_{z= 0     }^{\infty}
\left(\matrix{a_{ik} + zu^2 b_{ik} \cr z^{-2}b_{ik}}\right)z^ku^i $
in which every term is holomorphic in the $V$ chart except for 
$a_{1,-1}z .$
Subtracting the
holomorphic terms we are left with $T\sigma \sim \left(\matrix{ a_{1,-1}z \cr 0}\right)$
where $a_{1,-1} \in {\bf C}.$ 
Therefore $\sigma = T^{-1}T\sigma \sim  \left(\matrix{ a_{1,-1}z^{-1} \cr 0}\right)$
and we conclude that 
 $H^1(\ell_2, \widetilde{V}|{\ell_2}) $ is generated by 
$  \left(\matrix{ z^{-1} \cr 0}\right)$
and $l(R^1\pi_*\widetilde{V}) =  1.$

\subsection{General calculation proceedure}
Summing up, the numerical invariants
$l(Q)$ and $l(R^1\pi_*\widetilde{V})$  are calculated 
from the zero-th and first cohomologies 
of $\widetilde{V}$ on formal neighborhoods of the exceptional divisor.
Since bundles on $\widetilde{{\bf C}^2}$
with splitting type $j$ are determined by their restriction 
to the $(2j-2)-nd$
formal neighborhood, it turns out that the
invariants are determined by the cohomology groups 
 $H^0(\ell_{2j-2}, \widetilde{V}|{\ell_{2j-2}}) $
and $H^1(\ell_{2j-2}, \widetilde{V}|{\ell_{2j-2}}) $
which are in fact quite simple to calculate.

{\em Steps to calculate $l(Q)$ :}
 find the generators $\{\alpha_i\}$
of  $H^0(\ell_{2j-2}, \widetilde{V}|{\ell_{2j-2}}) ,$
determine the ${\cal O}_x$-module $M$ generated  
by the $\{\alpha_i\},$ write the natural inclusion 
$\rho: M \hookrightarrow M^{\vee\vee},$ 
then $l(Q)= dim \, coker\, \rho.$ 

{\em Steps to calculate  $l(R^1\pi_*\widetilde{V}):$}
 find the dimension of 
 $H^1(\ell_{2j-2}, \widetilde{V}|{\ell_{2j-2}}) $ as a 
$k(x)$-vector space, this dimension is
$l(R^1\pi_*\widetilde{V}).$

$\begin{array}{ll}
 E. Ballico &   E. Gasparim\\
Dept.\, of \, Mathematics & Depto.\, de \, Matematica\\
 University \, of \, Trento & Univ.\, Federal\,  de\,  Pernambuco\\
38050 Povo (TN) -Italy &  50670/901 Recife (PE)- Brasil \\
ballico@science.unitn.it &  gasparim@dmat.ufpe.br\hspace{1.2cm}\\
\end{array}$
\end{document}